\documentclass{amsart}
\usepackage{graphicx}
\usepackage{amscd}
\usepackage{amsmath}
\usepackage{amsfonts}
\usepackage{amssymb}
\theoremstyle{plain}
\newtheorem{theorem}{Theorem}[section]
\newtheorem{corollary}[theorem]{Corollary}
\newtheorem{lemma}[theorem]{Lemma}

\newtheorem{proposition}[theorem]{Proposition}
\theoremstyle{definition}
\newtheorem{example}[theorem]{Example}
\newtheorem{definition}[theorem]{Definition}

\newtheorem{remark}[theorem]{Remark}

\numberwithin{equation}{section}

\begin{document}
\title{$k$-hyponormality of powers of weighted shifts\\via Schur products}
\author{Ra{\'{u}}l Curto}
\address{Department of Mathematics\\
The University of Iowa\\
Iowa City, Iowa 52242\\
USA}
\email{curto@math.uiowa.edu}
\author{Sang Soo Park}
\address{Department of Mathematics\\
Kyungpook National University\\
Deagu, Korea }
\email{pss4855@hanmail.net}
\thanks{The research of the first named author was partially supported by NSF grant
DMS-9800931.\newline The research of the second named author was partially
supported by KOSEF research project no. R01-2000-00003.}

\begin{abstract}
We characterize $k$-hyponormality and quadratic hyponormality of powers of
weighted shifts using Schur product techniques.
\end{abstract}\maketitle

\section{Introduction}

Let $\mathcal{H}$ be a separable, infinite dimensional complex Hilbert space
and let $\mathcal{B(H)}$ be\ the algebra of bounded linear operators on
$\mathcal{H}$. An operator $T{\in\mathcal{B(H)}}$ is said to be
\textit{normal} if $T^{\ast}T=TT^{\ast}$, \textit{subnormal} if $T$ is the
restriction of a normal operator (acting on a Hilbert space $\mathcal{K}%
\supseteq\mathcal{H}$) to an invariant subspace, and \textit{hyponormal} if
$T^{\ast}T\geq TT^{\ast}$.

The Bram-Halmos criterion for subnormality states that an operator is
\textit{subnormal} if and only if%
\[
\sum\limits_{i,j}(T^{i}x_{j},T^{j}x_{i})\geq0
\]
for all finite collections $x_{0},x_{1},x_{2,}\cdots,x_{k}\in\mathcal{H}%
$\ (\cite{Bra}, \cite{Con}). Using Choleski's Algorithm for operator matrices,
it is easy to see that this is equivalent to the following positivity test:%
\begin{equation}
\left(
\begin{array}
[c]{llll}%
I & T^{\ast} & \cdots &  T^{\ast k}\\
T & T^{\ast}T & \cdots &  T^{\ast k}T\\
\vdots & \vdots & \ddots & \vdots\\
T^{k} & T^{\ast}T^{k} & \cdots &  T^{\ast k}T^{k}%
\end{array}
\right)  \geq0\;\;\,(\text{all }k\geq1). \label{eq1.1}%
\end{equation}
Condition (\ref{eq1.1}) provides a measure of the gap between hyponormality
and subnormality. \ The notion of $k$-hyponormality has been introduced and
studied in an attempt to bridge that gap (\cite{Ath}, \cite{BEJ}, \cite{Cu2},
\cite{CMX}, \cite{JL}, \cite{McCP}). In fact, the positivity condition
(\ref{eq1.1}) for $k=1$ is equivalent to the hyponormality of $T$, while
subnormality requires the validity of (\ref{eq1.1}) for all $k$. \ 

If we denote by $[A,B]:=AB-BA$\ the commutator of two operators $A\;$and $B,$
and if we define $T$ to be $k$\textit{-hyponormal} whenever the $k\times k$
operator matrix $M_{k}(T):=([T^{\ast j},T^{i}])_{i,j=1}^{k}$ is positive, or
equivalently, the $(k+1)\times(k+1)$operator matrix (\ref{eq1.1}) is positive,
then the Bram-Halmos criterion can be rephrased as saying that $T$\ is
subnormal if and only if $T$\ is $k$-hyponormal for every $k\geq1$ (\cite{CMX}).

Given a bounded sequence of positive numbers (called weights) $\alpha
\;:\;\alpha_{0},\alpha_{1},\alpha_{2},$ $\alpha_{3},\cdots$, the (unilateral)
weighted shift $W_{\alpha}$ associated with $\alpha$ is the operator on
$l^{2}(\mathbf{Z}_{+})$ defined by $W_{\alpha}e_{n}:=\alpha_{n}e_{n+1}$ for
all $n\geq0,$ where $\{e_{n}\}_{n=0}^{\infty}$ is the canonical orthonormal
basis for $l^{2}(\mathbf{Z}_{+}).$ It is straightforward to check that
$W_{\alpha}$ can never be normal, and that it is hyponormal if and only if
$\alpha_{n}\leq\alpha_{n+1}$ for all $n\geq0.$ The \textit{moments}
of\ $\alpha$ are usually defined by $\beta_{0}:=1,\beta_{n+1}:=\alpha_{n}%
\beta_{n}$\ ($n\geq0$) (\cite{Shi}); however, we will reserve this term for
the sequence $\gamma_{n}:=\beta_{n}^{2}\;(n\geq0)$. Berger's Theorem, which
follows, states that $W_{\alpha}$ is subnormal if and only if the moments of
$\alpha$ are the moments of a probability measure on $[0,\left\|  W_{\alpha
}\right\|  ^{2}]$.

\begin{theorem}
\label{berger}\textbf{(Berger's Theorem \cite{Con})} $W_{\alpha}$ is subnormal
if and only if there exists a Borel probability measure $\mu$ supported in
$\left[  0,{{||}W_{\alpha}{||}}^{2}\right]  $, with $||W_{\alpha}||^{2}%
\in\operatorname*{supp}\mu$ and such that%
\[
\gamma_{n}=\int t^{n}d\mu(t)\;\;(\text{all }n\geq0).
\]
\end{theorem}

In terms of $k$-hyponormality for weighted shifts, we will often use the
following basic result.

\begin{lemma}
\label{khyp} (\cite[Theorem 4]{Cu1}) $W_{\alpha}$ is $k$-hyponormal if and
only if the $(k+1)\times(k+1)$ Hankel matrices
\begin{equation}
A_{n,k}(\alpha):=\left(
\begin{array}
[c]{llll}%
\gamma_{n} & \gamma_{n+1} & \cdots & \gamma_{n+k}\\
\gamma_{n+1} & \gamma_{n+2} & \cdots & \gamma_{n+k+1}\\
\vdots & \vdots & \ddots & \vdots\\
\gamma_{n+k} & \gamma_{n+k+1} & \cdots & \gamma_{n+2k}%
\end{array}
\right)  \;(n\geq0) \label{eq1.2}%
\end{equation}
are all nonnegative.
\end{lemma}

In this article we study $k$-hyponormality and quadratic hyponormality of
powers of weighted shifts, using Schur product techniques. We characterize the
$k$-hyponormality of powers of $W_{\alpha}$ in terms of the $k$-hyponormality
of a finite collection of weighted shifts whose weight sequences are naturally
derived from $\alpha$. Similar techniques, when combined with the results in
\cite{BEJ}, \cite{Cu1}, \cite{RGWSI}, \cite{JP1} and \cite{JP2}, allow us to
deal with back-step extensions of weighted shifts, and with weak
$k$-hyponormality, including quadratic hyponormality.

\section{$k$-hyponormality of powers of weighted shifts}

For matrices $A,B\in$$M_{n}(\mathbf{C})$, we let $A\ast B$ denote their Schur
product. The followings result is well known.

\begin{lemma}
\label{lem2.1}\textbf{(\cite{Pau}) }If $A\geq0$ and $B\geq0$, then $A\ast
B\geq0.$

\begin{definition}
Let $\alpha\equiv\{\alpha_{n}\}_{n=0}^{\infty}$ and $\beta\equiv\{\beta
_{n}\}_{n=0}^{\infty}$. The Schur product of $\alpha$ and $\beta$ is defined
by $\alpha\beta:=\{{\alpha_{n}}{\beta_{n}}\}_{n=0}^{\infty}$.

\begin{theorem}
Let $\alpha\equiv\{\alpha_{n}\}_{n=0}^{\infty}$ and $\beta\equiv\{\beta
_{n}\}_{n=0}^{\infty}$ be two weight sequences, and assume that both
$W_{\alpha}$ and $W_{\beta}$ are\ $k$-hyponormal. Then $W_{\alpha\beta}$ is
$k$-hyponormal.

\begin{proof}
Let $\{\epsilon_{n}\}$ and $\{\eta_{n}\}$ be the moments of $\alpha$ and
$\beta$, respectively. By hypothesis, $A_{n,k}(\alpha)\geq0~$and
$A_{n,k}(\beta)\geq0$ (all $n\geq0$). Since the corresponding moments
$\gamma_{n}$ of $\alpha\beta$ satisfy $\gamma_{n}=\epsilon_{n}\eta_{n}~($all
$n\geq0)$, it follows that $A_{n,k}(\alpha\beta)=A_{n,k}(\alpha)\ast
A_{n,k}(\beta)\;\;$(all $n\geq0$). \ By Lemma\ \ref{lem2.1}, $A_{n,k}%
(\alpha\beta)\geq0$ (all $n\geq0$), so Lemma \ref{khyp} now implies that
$W_{\alpha\beta}$ is $k$-hyponormal.
\end{proof}
\end{theorem}
\end{definition}
\end{lemma}

\begin{corollary}
Let $W_{\alpha}$ and $W_{\beta}$ be two weighted shifts, and assume that each
is subnormal. \ Then $W_{\alpha\beta}$ is also subnormal.
\end{corollary}

\begin{proof}
This is a straightforward application of the Bram-Halmos Criterion.
\end{proof}

\begin{definition}
\label{defli}Given integers $i$ and $\ell$, with $\ell\geq1$ and $0\leq
i\leq\ell-1$, consider the decomposition $\mathcal{H}$ $\equiv l^{2}%
(\mathbf{Z}_{+})=\bigoplus_{j=0}^{\infty}\{e_{j}\}$, and define $\mathcal{H}%
_{i}:=$ $\bigoplus_{j=0}^{\infty}{e_{\ell j+i}}$. Moreover, for a weight
sequence $\alpha$ let $\alpha(\ell:i):=\{$ $\prod\limits_{m=0}^{\ell-1}%
\alpha_{\ell j+i+m}\}_{j=0}^{\infty}$. $\alpha(\ell:i)$ is the sequence of
products of weights in adjacent packets of size $\ell$, beginning with
$\alpha_{i}\cdot...\cdot\alpha_{i+\ell-1}$.

\begin{example}
Let $\alpha\equiv\{\alpha_{n}\}_{n=0}^{\infty}$ be a weight sequence. \ Then
\newline \textrm{(i)} $\ \alpha(2:0):\alpha_{0}\alpha_{1},\alpha_{2}\alpha
_{3},\alpha_{4}\alpha_{5},...$\newline \textrm{(ii) \ }$\alpha(3:1):\alpha
_{1}\alpha_{2}\alpha_{3},\alpha_{4}\alpha_{5}\alpha_{6},\alpha_{7}\alpha
_{8}\alpha_{9},...$\newline \textrm{(iii) \ }$\alpha(3:2):$ $\alpha_{2}%
\alpha_{3}\alpha_{4},\alpha_{5}\alpha_{6}\alpha_{7},\alpha_{8}\alpha_{9}%
\alpha_{10},...$

\begin{proposition}
\label{prop2.7}Let $\ell\geq1$, let $0\leq i\leq\ell-1$, and let $\alpha
(\ell:i)$ be as in Definition \ref{defli}. Then $W_{\alpha(\ell:i)}$ is
unitarily equivalent to $W_{\alpha}^{\ell}|_{\mathcal{H}_{i}}$. Therefore,
$W\,_{\alpha}^{\ell}$ is unitarily equivalent to $\bigoplus\limits_{i=0}%
^{\ell-1}W_{\alpha(\ell:i)}$.

\begin{proof}
\noindent Since $W_{\alpha}^{\ell}e_{\ell j+i}=\prod\limits_{m=0}^{\ell
-1}\alpha_{\ell j+i+m}{e}_{\ell(j+1)+i}$, it is clear that ${\mathcal{H}_{i}}$
is an invariant subspace for $W_{\alpha}^{\ell}$. Moreover, $(W_{\alpha}%
^{\ell})^{\ast}e_{\ell j+i}=\prod\limits_{m=0}^{\ell-1}\alpha_{\ell
(j-1)+i+m}{e}_{\ell(j-1)+i}$, so ${\mathcal{H}_{i}}$ is also invariant under
$(W_{\alpha}^{\ell})^{\ast}$. It follows that $\mathcal{H}_{i}$ is a reducing
subspace for $W_{\alpha}^{\ell}$. If we now define a unitary operator
$U:\mathcal{H}\longrightarrow\mathcal{H}_{i}$ by $U({e}_{j})$=$e_{\ell j+i}$,
we see at once that $U^{\ast}(W_{\alpha}^{\ell}|_{\mathcal{H}_{i}}%
)U$=$W_{\alpha(\ell:i)}$, as desired.
\end{proof}
\end{proposition}
\end{example}
\end{definition}

\begin{corollary}
\label{cor2.9}(a) $W_{\alpha}^{\ell}\;$is$\;k$-hyponormal $\Leftrightarrow$
$\;W_{\alpha(\ell:i)}$ is $k$-hyponormal for $0\leq i\leq\ell-1$.\newline (b)
$W_{\alpha}^{\ell}\;$is subnormal$\;\Leftrightarrow\;W_{\alpha(\ell:i)}$ is
subnormal for $0\leq i\leq\ell-1$.
\end{corollary}

Throughout the rest of this section, we assume that $W_{\alpha}$ is subnormal,
with Berger measure $\mu$. Observe that we can always write $\mu\equiv\nu
+\rho\delta_{0}$ where $\nu(\{0\})=0$, and that $W_{\alpha}^{\ell}$ is
subnormal whenever $W_{\alpha}$ is subnormal. By Corollary \ref{cor2.9}, we
know that each $W_{\alpha(\ell:i)}$ is subnormal, for $0\leq i\leq\ell-1$. We
now seek to identify the Berger measures $\mu_{i}$ corresponding to each
$W_{\alpha(\ell:i)}$.

\begin{theorem}
\label{thm2.10}(a) $d\mu_{0}(t)=d\mu(t^{1/\ell})$.\newline (b) For $1\leq
i\leq\ell-1$, $d\mu_{i}(t)=\frac{t^{i/\ell}}{\gamma_{i}}d\nu(t^{1/\ell})$.

\begin{proof}
Let ${\gamma_{n}}$ be the moments of $\alpha$ $(n\geq0)$. Then
\[
\int t^{n}d\mu_{0}(t)=\gamma_{\ell n}=\int t^{\ell n}d\mu(t),
\]
so $d{\mu}_{0}(t)=d\mu(t^{1/\ell})$. Similarly, for $1\leq i\leq l-1$,
\[
\int t^{n}d\mu_{i}(t)=\frac{\gamma_{\ell n+i}}{\gamma_{i}}=\frac{1}{\gamma
_{i}}\int t^{\ell n+i}d\nu(t),
\]
so $d\mu_{i}(t)=\frac{t^{i/\ell}}{\gamma_{i}}d\nu(t^{1/\ell})$.
\end{proof}
\end{theorem}

\section{Back-step Extensions of Weighted Shifts}

For a weight sequence $\alpha$, we consider the \textit{back-step extension}
$\alpha(x):x,\alpha_{0},\alpha_{1},$ $\alpha_{2},\alpha_{3},\cdots$ where
$x>0$.

\begin{lemma}
\label{lem2.11}Let $W_{\alpha}$ be a subnormal weighted shift with associated
Berger measure $\mu$.\newline (a) (cf. \cite[Proposition 8]{Cu1})
$W_{\alpha(x)}$ is subnormal if and only if (i) $\frac{1}{t}\in L^{1}(\mu)$
and (ii) $x^{2}\leq(||\frac{1}{t}||_{L^{1}(\mu)})^{-1}$. In particular,
$W_{\alpha(x)}$ is never subnormal when $\mu(\{0\})>0$.\newline (b) if
$x<(||\frac{1}{t}||_{L^{1}(\mu)})^{-1/2}$, the corresponding measure $\mu_{x}$
of $W_{\alpha(x)}$ satisfies $\mu_{x}(\{0\})$\newline $>0$. In particular,
$T:=W_{\alpha(\left\|  \frac{1}{t}\right\|  _{L^{1}(\mu)}^{-1/2})}$ is the
unique back-step extension of $W_{\alpha}$ with no mass at the origin.

\begin{proof}
(b) Let $\gamma_{n}$ be the moments of $T$. Since $T$ is subnormal, there
exists a Berger measure $\nu$ such that%
\[
\gamma_{n}=\int t^{n}d\nu=\left\{
\begin{array}
[c]{cc}%
1 & n=0\\
\frac{\int t^{n-1}d\mu}{\int\frac{1}{t}d\mu} & n\geq1
\end{array}
\right.  .
\]
Assume $x<(||\frac{1}{t}||_{L^{1}(\mu)})^{-1/2}$ and write $x=(1-\epsilon
)(||\frac{1}{t}||_{L^{1}(\mu)})^{-1/2}$ for $0<\epsilon<1$. The moments
$\eta_{n}$ of $W_{\alpha(x)}$ are such that%
\[
\eta_{n}=\left\{
\begin{array}
[c]{cc}%
1 & n=0\\
(1-\varepsilon)\gamma_{n}=\int t^{n}(1-\varepsilon)d\nu &  n\geq1
\end{array}
\right.  .
\]
It follows that $\mu_{x}=(1-\epsilon)d\nu+\epsilon\delta_{0}$.
\end{proof}
\end{lemma}

\begin{lemma}
\label{lem2.12}Let $W_{\alpha}$ be a subnormal weighted shift, let $\ell\geq
1$, and let $k\geq1$. The following statements are equivalent. \newline (a)
$W_{\alpha(x)}^{\ell}$ is $k$-hyponormal. \newline (b) $W_{\alpha(x)(\ell:0)}$
is $k$-hyponormal.

\begin{proof}
(a) $\Rightarrow$ (b). Straightforward from Corollary \ref{cor2.9}.
\newline (b) $\Rightarrow$ (a). By Corollary \ref{cor2.9}(b), we know that
$W_{\alpha(x)(\ell:i)}\equiv W_{a(\ell:i-1)}$ is subnormal, and by Proposition
\ref{prop2.7}, $W_{\alpha(x)}^{\ell}\cong\bigoplus\limits_{i=0}^{\ell
-1}W_{\alpha(x)(\ell:i)}$. It follows that for $1\leq i\leq l-1$,
$W_{\alpha(x)(\ell:i)}$ is $k$-hyponormal, which together with the assumption
that $W_{\alpha(x)(\ell:0)}$ is $k$-hyponormal shows that $W_{\alpha(x)}%
^{\ell}$ is $k$-hyponormal.
\end{proof}
\end{lemma}

\begin{theorem}
\label{thm2.13}Let $W_{\alpha}$ be subnormal, with Berger measure $\mu
\equiv\nu+\rho\delta_{0}$, and let $\ell\geq1$. Then $W_{\alpha(x)}^{\ell}$ is
subnormal if and only if $x\leq(||\frac{1}{t}||_{L^{1}(\nu)})^{-1/2}$. In
particular, if $\rho=0$, $W_{\alpha(x)}^{\ell}$ is subnormal if and only if
$W_{\alpha(x)}$ is subnormal.

\begin{proof}
It suffices to consider $W_{\alpha(x)(\ell:0)}$. Observe that $W_{\alpha
(x)(\ell:0)}$ is a back-step extension of $W_{\alpha(\ell:\ell-1)}$. By Lemma
\ref{lem2.11}, $W_{\alpha(x)(\ell:0)}$ is subnormal if and only if
$x^{2}\gamma_{\ell-1}\leq(||\frac{1}{t}||_{L^{1}(\mu_{\ell-1})})^{-1}%
=\gamma_{\ell-1}(|||\frac{1}{t}||_{L^{1}(\nu)})^{-1}$. Therefore
$W_{\alpha(x)(\ell:0)}$ is subnormal if and only if $x\leq(||{\frac{1}{t}%
}||_{L^{1}(\nu)})^{-1/2}$, as desired.
\end{proof}
\end{theorem}

\begin{remark}
Although for an operator $T$ the subnormality of $T^{\ell}$ does not imply the
subnormality of $T$, Theorem \ref{thm2.13} shows that this is the case for
back-step extensions of subnormal weighted shifts with Berger measures having
no mass at the origin.

\begin{theorem}
Let $W_{\alpha}$ be a subnormal weighted shift, with Berger measure $\mu$.
Then $W_{\alpha(x_{n},x_{n-1},\cdots,x_{1})}$ is subnormal if and only if
\newline (a) $\frac{1}{t^{j}}\in L^{1}(\mu)$ for all $1\leq j\leq n$,
\newline (b) $x_{1}\cdots x_{j}=(||\frac{1}{t^{j}}||_{L^{1}(\mu)})^{-1/2}$ for
$1\leq j\leq n-1$ and $x_{1}\cdots x_{n}\leq(||\frac{1}{t^{n}}||_{L^{1}(\mu
)})^{-1/2}$.

\begin{proof}
The case $n=1$ was established in \cite[Proposition 8]{Cu1}. Here, and without
loss of generality, we will only consider the case $n=2$.\newline
$(\Rightarrow)$ Assume that $W_{\alpha(x_{2},x_{1})}$ is subnormal. Since
$W_{\alpha(x_{1})}$ is a subnormal weighted shift possessing a subnormal
extension (namely $W_{\alpha(x_{2},x_{1})}$), Lemma \ref{lem2.11} implies that
$x_{1}=(||\frac{1}{t}||_{L^{1}(\mu)})^{-1/2}$. Moreover, since $W_{\alpha
(x_{2},x_{1})}$ is subnormal, we must have $W_{\alpha(x_{2},x_{1})}^{2}$
subnormal, so Lemma \ref{lem2.12} implies that $W_{\alpha(x_{2},x_{1}%
)(2:0)}\equiv W_{\alpha(2:0)(x_{2}x_{1})}$ is subnormal and $x_{1}x_{2}%
\leq\ (||\frac{1}{t}||_{L^{1}(\mu(t^{\frac{1}{2}}))})^{-1/2}=(||\frac{1}%
{t^{2}}||_{L^{1}(\mu)})^{-1/2}$. \newline $(\Leftarrow)$ Assume that (a) and
(b) hold. Since $\frac{1}{t}\in L^{1}(\mu)$ and $x_{1}^{2}=(||\frac{1}%
{t}||_{L^{1}(\mu)})$, we know that $W_{\alpha(x_{1})}$ is subnormal with
measure $\nu$ such that $\nu(\{0\})=0$. To check the subnormality of
$W_{\alpha(x_{2},x_{1})}=W_{\alpha(x_{1})(x_{2})}$, by Theorem \ref{thm2.13}
it suffices to check the subnormality of $W_{\alpha(x_{2},x_{1})}^{2}$, and by
Lemma \ref{lem2.12}, this reduces to verifying the subnormality of
$W_{\alpha(x_{2},x_{1})(2:0)}=W_{\alpha(2:0)(x_{2}x_{1})}$. If $\mu_{1}$
denotes the Berger measure of $W_{\alpha(2:0)}$, that is, $d\mu_{1}(t)\equiv
d\mu(t^{\frac{1}{2}})$, we know that $x_{2}x_{1}\leq(||\frac{1}{t^{2}%
}||_{L^{1}(\mu)})^{-1/2}=(||\frac{1}{t}||_{L^{1}(\mu_{1})})^{-1/2}$.
Therefore, we see that $W_{\alpha(2:0)(x_{2}x_{1})}$ is subnormal, using Lemma
\ref{lem2.11}. Thus, \newline $W_{\alpha(x_{2},x_{1})(2:0)}$ is subnormal, as desired.
\end{proof}
\end{theorem}
\end{remark}

\section{Some Revealing Examples}

Let $\mathcal{\alpha}:\alpha_{0},\alpha_{1},\alpha_{2},\cdots,\alpha
_{n},\cdots$ be a sequence of weights and let $\gamma_{n}$ be the
corresponding moments. For $x>0$ let ${\alpha}(x):x,{\alpha}_{0},{\alpha}%
_{1},\cdots$ be the associated back-step extension of $\alpha$ and assume that
$W_{\alpha}$ is subnormal. It follows from \cite[Theorem 4]{Cu1} that
$W_{\alpha(x)}$ is $k$-hyponormal if and only if
\[
D_{k}:=\left(
\begin{array}
[c]{ccccc}%
\frac{1}{x^{2}} & \gamma_{0} & \gamma_{1} & \cdots & \gamma_{k-1}\\
\gamma_{0} & \gamma_{1} & \gamma_{2} & \cdots & \gamma_{k}\\
\gamma_{1} & \gamma_{2} & \gamma_{3} & \cdots & \gamma_{k+1}\\
\vdots & \vdots & \vdots & \ddots & \vdots\\
\gamma_{k-1} & \gamma_{k} & \gamma_{k+1} & \cdots & \gamma_{2k-1}%
\end{array}
\right)  \geq0.
\]

\begin{theorem}
For $\ell\geq1$, $W_{\alpha(x)}^{\ell}$ is $k$-hyponormal if and only if
\[
D_{k;\ell}:=\left(
\begin{array}
[c]{ccccc}%
\frac{1}{x^{2}} & \gamma_{\ell-1} & \gamma_{2\ell-1} & \cdots & \gamma
_{k\ell-1}\\
\gamma_{\ell-1} & \gamma_{2\ell-1} & \gamma_{3\ell-1} & \cdots &
\gamma_{(k+1)\ell-1}\\
\gamma_{2\ell-1} & \gamma_{3\ell-1} & \gamma_{4\ell-1} & \cdots &
\gamma_{(k+2)\ell-1}\\
\vdots & \vdots & \vdots & \ddots & \vdots\\
\gamma_{k\ell-1} & \gamma_{(k+1)\ell-1} & \gamma_{(k+2)\ell-1} & \cdots &
\gamma_{2k\ell-1}%
\end{array}
\right)  \geq0.
\]

\begin{proof}
It suffices to check that $W_{\alpha(x)(\ell:0)}$ is $k$-hyponormal. Now, the
matrix detecting $k$-hyponormality for $W_{\alpha(x)(\ell:0)}$ is
\[
D_{k}=x^{2}\left(
\begin{array}
[c]{ccccc}%
\frac{1}{x^{2}} & \gamma_{\ell-1} & \gamma_{2\ell-1} & \cdots & \gamma
_{k\ell-1}\\
\gamma_{\ell-1} & \gamma_{2\ell-1} & \gamma_{3\ell-1} & \cdots &
\gamma_{(k+1)\ell-1}\\
\gamma_{2\ell-1} & \gamma_{3\ell-1} & \gamma_{4\ell-1} & \cdots &
\gamma_{(k+2)\ell-1}\\
\vdots & \vdots & \vdots & \ddots & \vdots\\
\gamma_{k\ell-1} & \gamma_{(k+1)\ell-1} & \gamma_{(k+2)\ell-1} & \cdots &
\gamma_{2k\ell-1}%
\end{array}
\right)  =x^{2}D_{k;\ell},\text{ }%
\]
so the result follows.
\end{proof}
\end{theorem}

\begin{proposition}
For $\ell\geq1$, let $\alpha:\sqrt{\frac{2}{3}},\sqrt{\frac{3}{4}},\sqrt
{\frac{4}{5}},\cdots$.\newline (1) $W_{\alpha(\sqrt{x})}^{\ell}$ is hyponormal
$\Leftrightarrow$ \ $x\leq\frac{(\ell+1)^{2}}{2(2\ell+1)}$ \newline (2)
$W_{\alpha(\sqrt{x})}^{\ell}$ is 2-hyponormal $\Leftrightarrow$ \ $x\;\leq
\frac{(\ell+1)^{2}(2\ell+1)^{2}}{2(3\ell+1)(4\ell^{2}+3\ell+1)}.$\newline (3)
$W_{\alpha(\sqrt{x})}^{\ell}$ is subnormal $\Leftrightarrow$ \ $x\;\leq
\frac{1}{2}.$
\end{proposition}

\begin{proof}
Observe that $\gamma_{k\ell-1}=\frac{2}{k\ell+1}$. Now consider
\[
D_{2}=\left(
\begin{array}
[c]{ccc}%
\frac{1}{x} & \frac{2}{\ell+1} & \frac{2}{2\ell+1}\\
\frac{2}{\ell+1} & \frac{2}{2\ell+1} & \frac{2}{3\ell+1}\\
\frac{2}{2\ell+1} & \frac{2}{3\ell+1} & \frac{2}{4\ell+1}%
\end{array}
\right)  .
\]
By direct calculation we obtain
\[
D_{2}\geq0\iff x\leq\frac{(\ell+1)^{2}(2\ell+1)^{2}}{2(3\ell+1)(4\ell
^{2}+3\ell+1)}.
\]
Moreover, since $W_{\alpha}$ is subnormal, with measure $2tdt$ (in particular,
with no mass at the origin), we see that $W_{\alpha(\sqrt{x})}^{\ell}$ is
subnormal $\iff$ $W_{\alpha(\sqrt{x})}$ is subnormal.
\end{proof}

\begin{corollary}
(\cite[Proposition 7]{Cu1}) Let $\alpha:\sqrt{\frac{2}{3}},\sqrt{\frac{3}{4}%
},\sqrt{\frac{4}{5}},\cdots$.\newline (1) $W_{\alpha(\sqrt{x})}$ is hyponormal
$\Leftrightarrow$ $x\;\leq\frac{2}{3}.$ \newline (2) $W_{\alpha(\sqrt{x})}$ is
2-hyponormal $\Leftrightarrow$ $x\;\leq\frac{9}{16}.$\newline (3)
$W_{\alpha(\sqrt{x})}$ is subnormal $\Leftrightarrow$ $x\;\leq\frac{1}{2}$.
\end{corollary}

\section{Quadratic hyponormality}

We recall some terminology and notation from \cite{Cu1}, \cite{RGWSI} and
\cite{RGWSII}. An operator $T$ is said to be \textit{quadratically hyponormal}
if $T+sT^{2}$ is hyponormal for every $s\in\mathbf{C}$. Let $W_{\alpha}$ be a
hyponormal weighted shift. For $s\in\mathbf{C}$, let $D(s):=[(W_{\alpha
}+sW_{\alpha}^{2})^{\ast},W_{\alpha}+sW_{\alpha}^{2}]$, let $P_{n}$ be the
orthogonal projection onto ${\bigvee}_{i=0}^{n}\{e_{i}\}$, and let%
\begin{align*}
D_{n}  &  \equiv D_{n}(s):=P_{n}[(W_{\alpha}+sW_{\alpha}^{2})^{\ast}%
,W_{\alpha}+sW_{\alpha}^{2}]P_{n}\\
&  =\left(
\begin{array}
[c]{cccccc}%
q_{0} & \bar{r}_{0} & 0 & \cdots & 0 & 0\\
r_{0} & q_{1} & \bar{r}_{1} & \cdots & 0 & 0\\
0 & r_{1} & q_{2} & \cdots & 0 & 0\\
\vdots & \vdots & \vdots & \ddots & \vdots & \vdots\\
0 & 0 & 0 & \cdots &  q_{n-1} & \bar{r}_{n-1}\\
0 & 0 & 0 & \cdots &  r_{n-1} & q_{n}%
\end{array}
\right)  ,
\end{align*}
where $q_{k}:=u_{k}+|s|^{2}v_{k},$ $r_{k}:=s\sqrt{w_{k}},$ $u_{k}:=\alpha
_{k}^{2}-\alpha_{k-1}^{2},$ $v_{k}:=\alpha_{k}^{2}\alpha_{k+1}^{2}%
-\alpha_{k-1}^{2}\alpha_{k-2}^{2},$ $w_{k}:=\alpha_{k}^{2}(\alpha_{k+1}%
^{2}-\alpha_{k-1}^{2})^{2}\ \ (k\geq0)$ and $\alpha_{-1}=\alpha_{-2}:=0.$
Clearly, $W_{\alpha}$ is quadratically hyponormal if and only if $D_{n}%
(s)\geq0$ for every $s\in\mathbf{C}$ and every $n\geq0$. Let $d_{n}%
(\cdot):=\det(D_{n}(\cdot))$. Then it follows from \cite{Cu1}, \cite{RGWSII}
that $d_{0}=q_{0}$, $d_{1}=q_{0}q_{1}-|r_{0}|^{2}$, and
\[
d_{n+2}=q_{n+2}d_{n+1}-|r_{n+1}|^{2}d_{n}\qquad(n\geq0),
\]
and that $d_{n}$ is actually a polynomial in $t:=|s|^{2}$ of degree $n+1$,
with Maclaurin expansion $d_{n}(t)\equiv\sum\limits_{i=0}^{n+1}c(n,i)t^{i}$.
This gives at once the relations%
\[%
\begin{array}
[c]{ccc}%
c(0,0)=u_{0} & c(0,1)=v_{0} & \\
c(1,0)=u_{1}u_{0} & c(1,1)=u_{1}v_{0}+u_{0}v_{1}-w_{0} & c(1,2)=v_{1}v_{0}%
\end{array}
,
\]
and%
\[
c(n+2,i)=u_{n+2}c(n+1,i)+v_{n+2}c(n+1,i-1)-w_{n+1}c(n,i-1)
\]%
\[
(n\geq0,0\leq i\leq n+1).
\]

To detect the positivity of $d_{n}$, the following notion was introduced in
\cite{RGWSII}.

\begin{definition}
We say that $W_{\alpha}$ is \textit{positively quadratically hyponormal} if
$c\left(  n,i\right)  \geq0$ $\,$for all $n,i\geq0$ with $0\leq i\leq n+1.$
\end{definition}

It is obvious that positive quadratic hyponormality implies quadratic
hyponormality; moreover, quadratic hyponormality does not necessarily imply
positive quadratic hyponormality \cite{JP1}.

\begin{proposition}
\label{prop4.2}With the above notation, assume that $u_{n+1}v_{n}\geq
w_{n}\;(n\geq3)$. Then $W_{\alpha}$ is positively quadratically hyponormal if
and only if $c(3,2)\geq0$ and $c(4,3)\geq0$.
\end{proposition}

\begin{proof}
Immediate from \cite[Corollary 3.3 and Theorem 3.9]{BEJ}.
\end{proof}

\begin{lemma}
Assume that $W_{\alpha}$ is subnormal and let $\ell\geq1$, $k\geq1$. The
following statements are equivalent. \newline (1) $W_{\alpha(x)}^{\ell}$ is
weakly $k$-hyponormal. \newline (2) $W_{\alpha(x)(\ell:0)}$ is weakly $k$-hyponormal.
\end{lemma}

\begin{proof}
Imitate the proof of Lemma \ref{lem2.12}.
\end{proof}

\begin{theorem}
Let $\alpha_{n}:=\sqrt{\frac{n+2}{n+3}}\;(n\geq0)$, let $\alpha\equiv
\{\alpha_{n}\}_{n=0}^{\infty}$, and let $\ell\geq1$. Then $W_{\alpha(\sqrt
{x})(\ell:0)}$ is positively quadratically hyponormal if and only if
\begin{equation}
\left\{
\begin{array}
[c]{cc}%
x\leq\frac{(\ell+1)^{2}}{2(2\ell+1)} & \ell=1,2\\
x\leq\frac{(\ell+1)^{2}(1+7\ell+34\ell^{2}+44\ell^{3})}{2(1+9\ell+45\ell
^{2}+99\ell^{3}+94\ell^{4})} & \ell\geq3
\end{array}
\right.  . \label{qhcriterion}%
\end{equation}
\end{theorem}

\begin{proof}
Let $\beta_{0}:=\sqrt{\frac{2x}{\ell+1}}$ and $\beta_{n}:=\sqrt{\frac{n\ell
+1}{(n+1)\ell+1}}~~(n\geq1)$. Then $W_{\alpha(\sqrt{x})(\ell:0)}=W_{\beta}$.
By direct calculation we see that
\[
u_{n}=\beta_{n}^{2}-\beta_{n-1}^{2}=\frac{\ell^{2}}{((n+1)\ell+1)(n\ell
+1)}\;\;(n\geq2),
\]
\newline
\[
v_{n}=\beta_{n}^{2}\beta_{n+1}^{2}-\beta_{n-1}^{2}\beta_{n-2}^{2}=\frac
{4\ell^{2}}{((n+2)\ell+1)(n\ell+1)}\;\;(n\geq3)
\]
and
\[
w_{n}=\beta_{n}^{2}(\beta_{n+1}^{2}-\beta_{n-1}^{2})^{2}=\frac{4\ell^{2}%
}{(n\ell+1)((n+1)\ell+1)((n+2)\ell+1)^{2}}\;\;(n\geq2).
\]
Since $W_{\beta}$ has the property $u_{n+1}v_{n}\geq w_{n}\;\;(n\geq3)$, by
Proposition \ref{prop4.2} it suffices to verify the nonnegativity of $c(3,2)$
and $c(4,3)$. By direct calculation,
\[
c(3,2)\geq0\iff\frac{(\ell+1)^{2}(7+11\ell)}{4(3+10\ell+11\ell^{2})}%
\]
and
\[
c(4,3)\geq0\iff\frac{(\ell+1)^{2}(1+7\ell+34\ell^{2}+44\ell^{3})}%
{2(1+9\ell+45\ell^{2}+99\ell^{3}+94\ell^{4})}.
\]
\newline On the other hand, the hyponormality condition for $W_{\beta}$ is
$x\leq\frac{(\ell+1)^{2}}{2(2\ell+1)}$. Finally, observe that
\[
\frac{(\ell+1)^{2}(7+11\ell)}{4(3+10\ell+11\ell^{2})}\geq\frac{(\ell+1)^{2}%
}{2(2\ell+1)}\;\;(\text{all }\ell\geq1),
\]%
\[
\frac{(\ell+1)^{2}(1+7\ell+34\ell^{2}+44\ell^{3})}{2(1+9\ell+45\ell^{2}%
+99\ell^{3}+94\ell^{4})}\geq\frac{(\ell+1)^{2}}{2(2\ell+1)}\;\;(\text{if }%
\ell=1,2),
\]
and
\[
\frac{(\ell+1)^{2}(1+7\ell+34\ell^{2}+44\ell^{3})}{2(1+9\ell+45\ell^{2}%
+99\ell^{3}+94\ell^{4})}\leq\frac{(\ell+1)^{2}}{2(2\ell+1)}\;\;(\text{if }%
\ell\geq3).
\]
This proves (\ref{qhcriterion}).
\end{proof}

\begin{corollary}
\label{cor5.5}Let $\alpha_{n}:=\sqrt{\frac{n+2}{n+3}}$ ($n\geq0$) and
$\alpha\equiv\{\alpha_{n}\}_{n=0}^{\infty}$. \newline (a) $W_{\alpha(\sqrt
{x})}^{2}$ is quadratically hyponormal $\iff$ $x\leq\frac{9}{10}.$\newline (b)
If $\ell\geq3$ and $x\leq\frac{(\ell+1)^{2}(1+7\ell+34\ell^{2}+44\ell^{3}%
)}{2(1+9\ell+45\ell^{2}+99\ell^{3}+94\ell^{4})}$, then $W_{\alpha(\sqrt{x}%
)}^{\ell}$ is quadratically hyponormal.

\begin{remark}
Let $\alpha$ be as in Corollary \ref{cor5.5}. \ Then $W_{\alpha(\sqrt{x})}%
^{2}$ is positively quadratically hyponormal if and only if $W_{\alpha
(\sqrt{x})}^{2}$ is quadratically hyponormal if and only if $x\leq\frac{9}%
{10}$. \ Moreover, for $x=\frac{9}{10}$, $W_{\alpha(\sqrt{x})}^{2}$ has the
first two weights equal, namely $\beta_{0}=\beta_{1}=\sqrt{\frac{3}{5}}$; this
example resembles \cite[Proposition 7]{Cu1}, where the first nontrivial
quadratically hyponormal weighted shift with two equal weights appears. (For
additional results along these lines, see \cite{CJ}.) \ Here we notice that
for $x=\frac{9}{10}$, not only $W_{\alpha(\sqrt{x})}$ is quadratically
hyponormal with two equal weights but also $W_{\alpha(\sqrt{x})}^{2}$ is
quadratically hyponormal! \ 
\end{remark}
\end{corollary}

\noindent\textbf{Acknowledgment.} \ The authors are indebted to the referee
for several helpful suggestions.


\begin{thebibliography}{McCP}
\bibitem[Ath]{Ath}A. Athavale, \emph{On joint hyponormality of operators,
}Proc. Amer. Math. Soc 103(1988), 417-42.

\bibitem[BEJ]{BEJ}J. Bae, G. Exner and I.B. Jung, \emph{Criteria for
positively quadratically hyponormal weighted shifts, }Proc. Amer. Math. Soc.,
to appear.

\bibitem[Bra]{Bra}J. Bram, \emph{Subnormal operators, }Duke Math. J 22(1955), 75-94.

\bibitem[Con]{Con}J.B. Conway, \emph{Subnormal Operators, }Research Notes in
Mathematics, vol. 51, Pitman Publ. Co., London, 1981.

\bibitem[Cu1]{Cu1}R. Curto, \emph{Quadratically hyponormal weighted shifts,
}Integral Equations Operator Theory 13(1990), 49-66.

\bibitem[Cu2]{Cu2}R. Curto, \emph{Joint hyponormality: A bridge between
hyponormality and subnormality, }Proc. Symposia Pure Math. 51(1990), 69-91.

\bibitem[CF1]{Houston}R. Curto and L. Fialkow, \emph{Recursiveness,
positivity, and truncated moment problems, }Houston J. Math. 17(1991), 603-635.

\bibitem[CF2]{RGWSI}R. Curto and L. Fialkow, \emph{Recursively generated
weighted shifts and the subnormal completion problem, }Integral Equations
Operator Theory 17(1993), 202-246.

\bibitem[CF3]{RGWSII}R. Curto and L. Fialkow, \emph{Recursively generated
weighted shifts and the subnormal completion problem, II, }Integral Equations
Operator Theory 18(1994), 369-426.

\bibitem[CJ]{CJ}R. Curto and I.B. Jung, \emph{Quadratically hyponormal
weighted shifts with two equal weights, }Integral Equations Operator Theory
37(2000), 208-231.

\bibitem[CJL]{CJL}R. Curto, I.B. Jung and W.Y. Lee, \emph{Extension and
extremality of recursively generated weighted shifts, }Proc. Amer. Math. Soc.,
to appear.

\bibitem[CP1]{CP1}R. Curto and M. Putinar, \emph{Existence of non-subnormal
polynomially hyponormal operators, }Bull. Amer. Math. Soc. 25(1991), 373-378.

\bibitem[CP2]{CP2}R. Curto and M. Putinar, \emph{Nearly subnormal operators
and moment problems, }J. Functional Analysis 115(1993), 480-497.

\bibitem[CMX]{CMX}R. Curto, P.S. Muhly and J. Xia, \emph{Hyponormal pairs of
commuting operators, }Operator Theory: Adv. Appl. 35(1988), 1-22.

\bibitem[JL]{JL}I.B. Jung and C. Li, \emph{A formula for $k$-hyponormality of
backstep extensions of subnormal weighted shifts, }Proc. Amer. Math. Soc., to appear.

\bibitem[JP1]{JP1}I.B. Jung and S. Park, \emph{Quadratically hyponormal
weighted shifts and their examples, }Integral Equations and Operator Theory
36(2000), 480-498.

\bibitem[JP2]{JP2}I.B. Jung and S. Park, \emph{Cubically hyponormal weighted
shifts and their examples, }J. Math. Anal. Appl. 247(2000), 557-569.

\bibitem[McCP]{McCP}S. McCullough and V. Paulsen, \emph{A note on joint
hyponormality, }Proc. Amer. Math. Soc. 107(1989), 187-195.

\bibitem[Pau]{Pau}V. Paulsen, \emph{Completely bounded maps and dilations,
}Pitman Research Notes in Mathematics Series, vol. 146, Longman Sci. Tech.,
New York, 1986.

\bibitem[Sta]{Sta}J. Stampfli, \emph{Which weighted shifts are subnormal,
}Pacific J. Math. 17(1996), 367-379.

\bibitem[Shi]{Shi}A. Shields, \emph{Weighed shift operators and analytic
function theory, }Proc. Amer. Math. Soc. 107(1989), 187-195.

\bibitem[Wol]{Wol}Wolfram Research Inc., \textit{Mathematica} Version 3.0,
Wolfram Research Inc., Champaign, IL, 1996.
\end{thebibliography}
\end{document}